\documentstyle[amscd,amssymb,verbatim,epsf]{amsart}

\begin{document}
\theoremstyle{plain}
\newtheorem{Thm}{Theorem}
\newtheorem{Cor}{Corollary}
\newtheorem{Ex}{Example}
\newtheorem{Con}{Conjecture}
\newtheorem{Main}{Main Theorem}
\newtheorem{Lem}{Lemma}
\newtheorem{Prop}{Proposition}

\theoremstyle{definition}
\newtheorem{Def}{Definition}
\newtheorem{Note}{Note}

\theoremstyle{remark}
\newtheorem{notation}{Notation}
\renewcommand{\thenotation}{}

\errorcontextlines=0
\numberwithin{equation}{section}
\renewcommand{\rm}{\normalshape}%

\title[Area-stationary surfaces]%
   {On area-stationary surfaces in certain neutral K\"ahler 4-manifolds}
\author{Brendan Guilfoyle}
\address{Brendan Guilfoyle\\
          Department of Mathematics and Computing \\
          Institute of Technology, Tralee \\
          Clash \\
          Tralee  \\
          Co. Kerry \\
          Ireland.}
\email{brendan.guilfoyle@@ittralee.ie}
\author{Wilhelm Klingenberg}
\address{Wilhelm Klingenberg\\
 Department of Mathematical Sciences\\
 University of Durham\\
 Durham DH1 3LE\\
 United Kingdom.}
\email{wilhelm.klingenberg@@durham.ac.uk }

\date{10th November, 2006}

\begin{abstract}
We study surfaces in TN that are area-stationary with respect to a neutral K\"ahler metric constructed on TN from a 
riemannian metric g on N. We show that holomorphic curves in TN are area-stationary, while lagrangian surfaces that
are area-stationary are also holomorphic and hence totally null. However, in general, area stationary surfaces are
not holomorphic. We prove this by constructing counter-examples. In the case where g is rotationally symmetric, we 
find all area stationary surfaces that arise as graphs of sections of the bundle TN$\rightarrow$N and that are 
rotationally symmetric. When (N,g) is the round 2-sphere, TN can be identified with the space of oriented affine lines
in ${\Bbb{R}}^3$, and we exhibit a two parameter family of area-stationary tori that are neither holomorphic nor 
lagrangian.
\end{abstract}

\maketitle

One-dimensional submanifolds of neutral K\"ahler four-manifolds have been studied recently in the context of twistor 
theory and integrable systems ({\it cf}. \cite{DunWest} and references therein). For example, quotienting out by the 
integral curves of non-null or null Killing vector fields leads to Einstein-Weyl three-manifolds or projective surfaces, 
respectively. 

In the case of two-dimensional submanifolds in neutral K\"ahler four-manifolds, the objects of study
in this paper, the situation is more complex. In particular, the metric induced on such a submanifold by the neutral 
metric can be positive or negative definite, Lorentz or degenerate, with two possible degrees of degeneracy. Moreover,
while the geometry of surfaces in positive definite K\"ahler four-manifolds is well developed, particularly in the 
K\"ahler-Einstein case \cite{ShoenWolf}, for indefinite metrics many of these results do not hold.

The main purpose of this paper then is two-fold: to investigate the geometric properties of two-dimensional 
submanifolds of a class of neutral K\"ahler four-manifolds and, by so doing, to illustrate the differences between 
the neutral and Hermitian cases. 

The particular class of neutral K\"ahler structures we consider have recently been studied on TN, the total space of 
the tangent bundle to a riemannian two-manifold (N,g) \cite{gak4} \cite{gak5}. This construction is motivated by the 
neutral K\"ahler metric on the space of oriented lines in ${\Bbb{R}}^3$
and on the space of time-like lines in ${\Bbb{R}}^{2+1}$. Aside from the signature, these K\"ahler four-manifolds differ 
from the more commonly studied
K\"ahler four-manifolds in a number of crucial ways: they are non-compact and are 
K\"ahler-Einstein only in the case when (N,g) is flat. They are, however, scalar flat, and the symplectic structure is 
exact. 

In the next section we discuss neutral K\"ahler metrics and the some of their properties. We also outline the
construction of the neutral K\"ahler structure on TN and the geometric structures induced on surfaces in TN. In the 
following section we derive the equations for area-stationary graphs in TN and show that holomorphic curves are
area-stationary. We also show that lagrangian area-stationary graphs are totally null. In section 3 we look at
rotationally symmetric graphs and determine all of these that area-stationary. In addition, we give a construction for
surfaces on which the induced metric is degenerate at every point.

In the final section, we look at the case of TS$^2$, which we identify with the space ${\Bbb{L}}$ of oriented affine 
lines in ${\Bbb{R}}^3$. There we construct area-stationary tori that are neither holomorphic nor lagrangian, and 
investigate their geometric properties. This two-parameter family of surfaces are analogous to the catenoid in 
${\Bbb{R}}^3$, being the 
unique rotationally symmetric area-stationary surfaces in ${\Bbb{L}}$.

\section{The Neutral K\"ahler Metric on TN}

We begin with some general properties of a K\"ahler surface (${\Bbb{M}},{\Bbb{G}},{\Bbb{J}},\Omega$). 
That is, ${\Bbb{M}}$ is a real 4-manifold 
endowed with the following structures. First, there is the metric ${\Bbb{G}}$, which we do not insist be positive definite
- it may also have neutral signature ($++--$). In order to deal with both cases simultaneously we assume that the metric 
can be diagonalised pointwise to ($1,1,\epsilon,\epsilon$), for $\epsilon=\pm1$.

In addition, we have a complex structure ${\Bbb{J}}$, which is a mapping
${\Bbb{J}}:\mbox{T}_p{\Bbb{M}}\rightarrow \mbox{T}_p{\Bbb{M}}$ at each $p\in{\Bbb{M}}$, which satisfies ${\Bbb{J}}^2=-{\mbox{Id}}$ and 
an integrability condition. 
Finally, there is a symplectic form $\Omega$, which is a closed, non-degenerate 2-form. These structures are required  
to satisfy the compatibility conditions:
\[
{\Bbb{G}}({\Bbb{J}}\cdot,{\Bbb{J}}\cdot)={\Bbb{G}}(\cdot,\cdot) \qquad\qquad
{\Bbb{G}}(\cdot,\cdot)=\Omega({\Bbb{J}}\cdot,\cdot).
\]

The following calibration identity highlights the difference between the case where ${\Bbb{G}}$ is positive definite and 
where it is neutral.

\begin{Thm}\cite{gak5}\label{t:wirt}
Let $p\in{\Bbb{M}}$ and $v_1,v_2\in T_p{\Bbb{M}}$ span a plane. Then
\[
\Omega(v_1,v_2)^2+\epsilon\varsigma^2(v_1,v_2)={\mbox{det }}{\Bbb{G}}(v_i,v_j),
\]
where $\varsigma^2(v_1,v_2)\geq0$ with equality iff $\{v_1,v_2\}$ spans a complex plane.
\end{Thm}

In the Hermitian case, $\epsilon=1$ and the above equality implies Wirtinger's inequality: the symplectic area is
bounded above by the metric area. 

Given a surface $\Sigma$ in ${\Bbb{M}}$, we say that $\Sigma$ is {\it holomorphic} if ${\Bbb{J}}$ preserves the tangent
space of $\Sigma$, while it is {\it lagrangian} if the symplectic form pulled back to $\Sigma$ vanishes.
A further consequence of the above Theorem is that, in the positive definite case, a surface cannot be both 
holomorphic and lagrangian. In the neutral case, however, this is not true: a surface can be both holomorphic and 
lagrangian, the only requirement being that the metric must be maximally degenerate along such a surface. We call such
surfaces {\it totally null} surfaces and the full details of this are included in Proposition \ref{p:indmet} below.

We turn now to the construction of a neutral K\"ahler structure on TN - further details can be found in \cite{gak4} 
\cite{gak5}. 
Given a riemannian 2-manifold (N,g,j) we construct
a canonical K\"ahler structure (${\Bbb{J}}$,$\Omega$,${\Bbb{G}}$) on the tangent bundle TN as follows. The
Levi-Civita connection associated with g splits the tangent bundle TTN$\cong$TN$\oplus$TN and the almost complex 
structure is defined to be ${\Bbb{J}}=\mbox{j}\oplus \mbox{j}$. This turns out to satisfy the appropriate 
integrability condition and so we have a complex structure on TN.

To define the symplectic form, consider the metric g as a mapping from TN to T$^*$N and pull back the canonical 
symplectic 2-form $\Omega^*$ on T$^*$N to a symplectic 2-form $\Omega$ on TN. Finally, the metric is defined as above by
${\Bbb{G}}(\cdot,\cdot)=\Omega({\Bbb{J}}\cdot,\cdot)$. The triple (${\Bbb{J}}$, $\Omega$, ${\Bbb{G}}$) 
determine a K\"ahler structure on TN. 

\begin{Prop}\cite{gak4}
Let (TN,${\Bbb{J}}$,$\Omega$,${\Bbb{G}}$) be the K\"ahler surface, as above. Then the metric ${\Bbb{G}}$ has
neutral signature ($++--$) and is scalar-flat. Moreover, ${\Bbb{G}}$ is K\"ahler-Einstein iff $g$ is flat, and
${\Bbb{G}}$ is conformally flat iff $g$ is of constant curvature.
\end{Prop}

Choose holomorphic coordinates $\xi$ on N so that
$ds^2=e^{2u}d\xi d\bar{\xi}$ for $u=u(\xi,\bar{\xi})$, and corresponding coordinates ($\xi$,$\eta$) on TN  by
identifying
\[
(\xi,\eta) \leftrightarrow \eta\frac{\partial}{\partial \xi}+\bar{\eta}\frac{\partial}{\partial \bar{\xi}}
                   \in \mbox{T}_\xi \mbox{N}.
\] 
These coordinates turn out to be holomorphic with respect to the above complex structure on TN and the symplectic 2-form 
has the following expression:

\[
\Omega=2{\Bbb{R}}\mbox{e}\left(e^{2u}d\eta\wedge d\bar{\xi}+\eta\partial(e^{2u})d\xi\wedge d\bar{\xi}\right).
\]
Here we have introduced the notation $\partial$ for differentiation with respect to $\xi$ - notation that we
will use throughout this paper.
The symplectic 2-form is globally exact $\Omega=d\Theta$, where $\Theta=\eta e^{2u}d\bar{\xi}+\bar{\eta}e^{2u}d\xi$.
Thus, for a closed surface $\Sigma$ in TN
\[
\int_\Sigma\Omega=0.
\]

On the other hand, the K\"ahler metric ${\Bbb{G}}$ is given in holomorphic coordinates by
\[
{\Bbb{G}}=2{\Bbb{I}}\mbox{m}\left(e^{2u}d\bar{\eta} d\xi-\eta\partial(e^{2u})d\xi d\bar{\xi}\right).
\]

We now consider surfaces in TN which arise as the graph of a local section of the bundle 
TN$\rightarrow$N,
that is, a map $\xi\rightarrow(\xi,\eta=F(\xi,\bar{\xi}))$. For such a surface introduce the following notation for
the complex slopes:
\[
\sigma=-\partial \bar{F} \qquad\qquad \rho=e^{-2u}\partial\left(Fe^{2u}\right),
\]
and let $\lambda={\Bbb{I}}{\mbox{m}}\;\rho$.

\begin{Prop}\cite{gak4}
A graph of a local section is holomorphic iff $\sigma=0$ and is lagrangian iff $\lambda=0$
\end{Prop}

Turning to the metric on a graph, the following makes explicit the identity contained in Theorem \ref{t:wirt}:

\begin{Prop}\label{p:indmet}
The metric (and its inverse) induced on the graph of a section by the K\"ahler metric is given in coordinates ($\xi,\bar{\xi}$) by;
\[
{\Bbb{G}}=e^{2u}\left[\begin{matrix}
i\sigma & -\lambda\\
-\lambda & -i\bar{\sigma}\\
\end{matrix}
\right]
\qquad\qquad
{\Bbb{G}}^{-1}=\frac{e^{-2u}}{\lambda^2-\sigma\bar{\sigma}}\left[\begin{matrix}
i\bar{\sigma} & -\lambda\\
-\lambda & -i\sigma\\
\end{matrix}
\right].
\]

In particular, the determinant of the induced metric is 
$|{\Bbb{G}}|=(\lambda^2-\sigma\bar{\sigma})e^{4u}$. Thus, the metric
is lorentz iff $\lambda^2<\sigma\bar{\sigma}$, riemannian iff $\lambda^2>\sigma\bar{\sigma}$ and degenerate if
$\lambda^2=\sigma\bar{\sigma}$. The metric induced on a holomorphic, lagrangian graph is identically zero, and
we call such a surface totally null.
\end{Prop}

This has the following corollary:

\begin{Cor}
The metric induced on a closed surface in TN cannot be positive (or negative) definite everywhere. 
\end{Cor}
\begin{pf}
Since the symplectic form is exact, as noted previously, its integral over any closed surface is zero. Thus, the
symplectic form must vanish somewhere on the surface. By Theorem \ref{t:wirt}, at such points the determinant of
the induced metric is either zero or negative. That is, the metric must be either degenerate or lorentz at these points.
\end{pf}

\section{Area-stationary Graphs}

The area form of the induced metric is  $|{\Bbb{G}}|^{\scriptstyle{\frac{1}{2}}}d\xi\wedge d\bar{\xi}$, and the 
following proposition
deals with stationary values of the area functional:

\begin{Prop}
A surface $\Sigma\hookrightarrow TN$ which is given by the graph of a function 
$\xi\rightarrow(\xi,\eta=F(\xi,\bar{\xi}))$ is area-stationary iff
\begin{equation}\label{e:max1}
i\partial\left(\frac{\lambda}{\sqrt{|\lambda^2-\sigma\bar{\sigma}|}}\right)
    -e^{-2u}\bar{\partial}\left(\frac{\sigma e^{2u}}{\sqrt{|\lambda^2-\sigma\bar{\sigma}|}}\right)=0.
\end{equation}
\end{Prop}
\begin{pf}
From Proposition \ref{p:indmet} the area functional evaluated on a graph $\Sigma$ is
\[
A(\Sigma)=\int_\Sigma d\mbox{vol}({\Bbb{G}})=\int_\Sigma |{\Bbb{G}}|^{\scriptstyle{\frac{1}{2}}}d\xi d\bar{\xi}
    =\int_\Sigma |\lambda^2-\sigma\bar{\sigma}|^{\scriptstyle{\frac{1}{2}}}e^{2u}\;d\xi d\bar{\xi}.
\] 
Varying the graph F we have
\[
\delta A(F)={\textstyle{\frac{1}{2}}}\int_\Sigma |\lambda^2-\sigma\bar{\sigma}|^{\scriptstyle{-\frac{1}{2}}}
         \delta(\lambda^2-\sigma\bar{\sigma})e^{2u}\;d\xi d\bar{\xi}.
\]
Now
\begin{align}
\delta(\lambda^2-\sigma\bar{\sigma})=&2\lambda\delta\lambda-\sigma\delta\bar{\sigma}-\bar{\sigma}\delta\sigma
   \nonumber\\
   & =-i\lambda e^{-2u}\left(\partial(\delta F e^{2u})-\bar{\partial}(\delta \bar{F} e^{2u})\right)
      +\sigma\bar{\partial}(\delta F)+\bar{\sigma}\partial(\delta\bar{F})\nonumber,
\end{align}
and so
\[
\delta A(F)={\textstyle{\frac{1}{2}}}\int_\Sigma 
         \left[-i\lambda e^{-2u}\left(\partial(\delta F e^{2u})-\bar{\partial}(\delta \bar{F} e^{2u})\right)
      +\sigma\bar{\partial}(\delta F)+\bar{\sigma}\partial(\delta\bar{F})\right]\frac{e^{2u}\;d\xi d\bar{\xi}}
     {|\lambda^2-\sigma\bar{\sigma}|^{\scriptstyle{\frac{1}{2}}}}.
\]
Integrating by parts we get
\begin{align}
\delta A(F)=&{\textstyle{\frac{1}{2}}}\int_\Sigma \left(ie^{2u}\partial\left(\frac{\lambda}
              {|\lambda^2-\sigma\bar{\sigma}|^{\scriptstyle{\frac{1}{2}}}}\right)
 -\bar{\partial}\left(\frac{\sigma e^{2u}}{|\lambda^2-\sigma\bar{\sigma}|^{\scriptstyle{\frac{1}{2}}}})\right)\right)\delta F\;d\xi d\bar{\xi}
\nonumber\\
&+{\textstyle{\frac{1}{2}}}\int_\Sigma \left(-ie^{2u}\bar{\partial}\left(\frac{\lambda}
              {|\lambda^2-\sigma\bar{\sigma}|^{\scriptstyle{\frac{1}{2}}}}\right)
 -\partial\left(\frac{\bar{\sigma}e^{2u}}{|\lambda^2-\sigma\bar{\sigma}|^{\scriptstyle{\frac{1}{2}}}})\right)\right)\delta \bar{F}\;d\xi d\bar{\xi} \nonumber.
\end{align}
A graph $F$ is area-stationary if $\delta A(F)=0$ for all $\delta F$, and so the result follows.
\end{pf}

The previous Proposition has the following Corollary:

\begin{Cor}
Holomorphic graphs in TN are area-stationary, while lagrangian graphs that
are area-stationary are also holomorphic and, hence, totally null.
\end{Cor}
\begin{pf}
Suppose that the graph of the section is holomorphic. Then $\sigma=0$ and we see that
\[
i\partial\left(\frac{\lambda}{\sqrt{|\lambda^2-\sigma\bar{\sigma}|}}\right)
    -e^{-2u}\bar{\partial}\left(\frac{\sigma e^{2u}}{\sqrt{|\lambda^2-\sigma\bar{\sigma}|}}\right)=
i\partial\left(\frac{\lambda}{|\lambda|}\right)=0,
\]
and so by the previous Proposition it is area-stationary.

On the other hand, a lagrangian graph has $\lambda=0$, and so, if, in addition, it is area-stationary
\[
0=i\partial\left(\frac{\lambda}{\sqrt{|\lambda^2-\sigma\bar{\sigma}|}}\right)
    -e^{-2u}\bar{\partial}\left(\frac{\sigma e^{2u}}{\sqrt{|\lambda^2-\sigma\bar{\sigma}|}}\right)=
-e^{-2u}\bar{\partial}\left(\frac{\sigma e^{2u}}{|\sigma|}\right),
\]
so that $\sigma e^{2u}|\sigma|^{-1}$ is holomorphic. This is impossible unless $\sigma=0$, in which case, the 
graph is holomorphic, as claimed.
\end{pf}

\section{Rotationally Symmetric Area Stationary Graphs}

Let (N,g) be a riemannian two-manifold.

\begin{Def}
The metric g is {\it rotationally symmetric} if there exists a conformal coordinate system ($\xi,\bar{\xi}$) such that
$g=e^{2u}d\xi d\bar{\xi}$ for $u=u(|\xi|)$.
\end{Def}

In other words, the metric is invariant under $\xi\rightarrow\xi e^{iC}$. Such an isometry of
(N,g) lifts to an isometry $(\xi,\eta)\rightarrow(\xi e^{iC},\eta e^{iC})$ of the K\"ahler metric on TN by the 
derivative map \cite{gak5}.
\begin{Def}
Let g be rotationally symmetric. A surface in TN is {\it rotationally symmetric} if it is invariant under the induced
isometry of TN.
\end{Def}

\begin{Lem}
A graph $\xi\rightarrow(\xi,\eta=F(\xi,\bar{\xi}))$ is rotationally symmetric iff $F(\xi,\bar{\xi})=G(R)e^{i\theta}$
for some complex-valued function $G$, where $\xi=Re^{i\theta}$.
\end{Lem}

The following Theorem characterises area-stationary graphs in TN that are rotationally symmetric:
 
\begin{Thm}
Let ($N,g$) be a rotationally symmetric riemannian two-manifold and ${\Bbb{G}}$ be the associated neutral K\"ahler metric 
on $TN$. A rotationally symmetric surface which is given by the graph of a local section 
$\xi\rightarrow (\xi,\eta=F(\xi,\bar{\xi}))$ is area-stationary with respect to ${\Bbb{G}}$ iff
\[
F=\left[A_1R+B_1R^{-1}e^{-2u}\pm i\left[A_2R^2+B_2e^{-2u}-B_1^2R^{-2}e^{-4u}\right]^{\scriptstyle{\frac{1}{2}}}\right]e^{i\theta},
\]
for $A_1,A_2,B_1,B_2\in{\Bbb{R}}$, $A_2\neq0$, where $\xi=Re^{i\theta}$ and $g=e^{2u}d\xi d\bar{\xi}$ for $u=u(R)$.
\end{Thm}
\begin{pf}
Let $F=(H\pm i\Psi^{\scriptstyle{\frac{1}{2}}})e^{i\theta}$ for real functions $H=H(R)$ and $\Psi=\Psi(R)$. 
Substituting this in equation (\ref{e:max1}) we get a pair of coupled non-linear 2nd order ordinary differential 
equations for $H$ and $\Psi$, which can be written
\begin{equation}\label{e:eq1&2}
\ddot{\Psi}+p_1\dot{\Psi}+q_1\Psi=L_1
\qquad\quad
\ddot{\Psi}+p_2\dot{\Psi}+q_2\Psi=L_2,
\end{equation}
where a dot represents differentiation with respect to $R$ and 
\[
p_1=-\frac{1+R^2(\ddot{u}-2\dot{u}^2)}{R(1+R\dot{u})}
\qquad\qquad
q_1=-\frac{2(\dot{u}-R(\ddot{u}-2\dot{u}^2))}{R(1+R\dot{u})}
\]
\[
L_1=\frac{R\dot{H}-H}{R^2(1+R\dot{u})^2}\left[R^2(1+R\dot{u})\ddot{H}-(1+2R\dot{u}+R^2\ddot{u})(R\dot{H}-H)\right],
\]
and
\[
p_2=-\frac{2R\ddot{H}}{R\dot{H}-H}-\frac{3+4R\dot{u}^2-R^2(\ddot{u}-2\dot{u}^2)}{R(1+R\dot{u})}
\]
\[
q_2=-\frac{4R\ddot{H}}{R\dot{H}-H}-\frac{2(3\dot{u}+R(6\dot{u}^2-\ddot{u})-2R^2(\ddot{u}-2\dot{u}^2)\dot{u})}
                  {R(1+R\dot{u})}
\]
\[
L_2=-\frac{2(R\dot{H}-H)^2}{R^2},
\]
To solve these equations proceed as follows: first solve the homogenous version of the first equation in 
(\ref{e:eq1&2}) for
$\Psi$ and then use variation of parameters to solve the inhomogeneous equation for 
$\Psi=\Psi(R,\dot{u},\ddot{u},H,\dot{H},\ddot{H})$. Then substitute this in the second equation of 
(\ref{e:eq1&2}) and solve for $H=H(R)$.

We now carry this out in detail. Start by noting that $\Psi=R^2$ is a solution of the homogenous version of the first equation of (\ref{e:eq1&2}). Now the other 
linearly independent solution of the homogenous equation can be found by recourse to the following lemma:

\begin{Lem}\cite{BirRota}\label{l:2ode}
Let $\Psi=\Psi_1$ be a solution of the 2nd order linear homogenous ordinary differential equation
$\ddot{\Psi}+p(R)\dot{\Psi}+q(R)\Psi=0$. Then the other linearly independent solution is
\[
\Psi_2=\Psi_1\int \Psi_1^{-2}e^{-P}dr \qquad \qquad \mbox{where} \qquad P(R)=\int p(R)dR.
\]
\end{Lem}
In our case, $\Psi_1=R^2$ and 
\[
p(R)=p_1=-\frac{1+R^2(\ddot{u}-2\dot{u}^2)}{R(1+R\dot{u})}=-\frac{d}{dR}\ln\left[R(1+R\dot{u})e^{-2u}\right],
\]
so that the second solution is
\[
\Psi_2=R^2\int R^{-3}(1+R\dot{u})e^{-2u}dR=-{\textstyle{\frac{1}{2}}}e^{-2u}.
\]
Thus the homogenous solution to first equation of (\ref{e:eq1&2}) is
\[
\Psi=A_2R^2+B_2e^{-2u},
\]
for real constants $A_2$ and $B_2$.

To solve the full equation we now use variation of parameters:
\[
\Psi=R^2(A_2-I_1)+e^{-2u}(B_2+I_2),
\]
where
\[
I_1=\int\frac{L_1}{2R(1+R\dot{u})}dR \qquad\qquad I_2=\int\frac{RL_1}{2(1+R\dot{u})e^{-2u}}dR  .
\]
The first of these can be completely integrated
\[
I_1=-\left[\frac{R\dot{H}-H}{2R(1+R\dot{u})}\right]^2,
\]
while the second can be integrated by parts to
\[
I_2=-\left[\frac{R\dot{H}-H}{2(1+R\dot{u})}\right]^2e^{2u}+\int\frac{(R\dot{H}-H)^2}{2R(1+R\dot{u})}e^{2u}dR.
\]
Thus the solution of the first equation of (\ref{e:eq1&2}) is 
\begin{equation}\label{e:psi}
\Psi=A_2R^2+B_2e^{-2u}+e^{-2u}\int\frac{(R\dot{H}-H)^2}{2R(1+R\dot{u})}e^{2u}dR.
\end{equation}

Substituting this in the second equation of (\ref{e:eq1&2}) yields the following:
\[
A_2(1+R\dot{u})\left[R^2(1+R\dot{u})\ddot{H}+[1+2R\dot{u}-R^2(\ddot{u}-2\dot{u}^2)]\left[R\dot{H}-H\right]  \right]=0.
\]
If $A_2=0$ we find that the surface is degenerate, ({\it cf}. Proposition \ref{p:deg} below). Also, since $1+R\dot{u}$ 
is not identically zero, we must solve
\[
R^2(1+R\dot{u})\ddot{H}+[1+2R\dot{u}-R^2(\ddot{u}-2\dot{u}^2)]\left[R\dot{H}-H\right]=0.
\]
This has one solution given by $H_1=R$ and we apply Lemma \ref{l:2ode} to find the second solution:
\[
H_2=R\int R^{-2}e^{-P}dR,
\]
where
\[
P=\int \frac{1+2R\dot{u}-R^2(\ddot{u}-2\dot{u}^2)}{R(1+R\dot{u})}dR=-\ln\left(R^{-1}(1+R\dot{u})e^{-2u}\right).
\]
Thus 
\[
H_2=R\int R^{-3}(1+R\dot{u})e^{-2u}dR=-{\textstyle{\frac{1}{2}}}R^{-1}e^{-2u},
\]
and the complete solution is
\[
H=A_1R+B_1R^{-1}e^{-2u}.
\]
Substituting this back in equation (\ref{e:psi}) we find that
\[
\Psi=A_2R^2+B_2e^{-2u}-B_1^2R^{-2}e^{-4u},
\]
which completes the theorem.
\end{pf}

The following deals with the case $A_2=0$:
\begin{Prop}\label{p:deg}
Let ($N,g$) be a rotationally symmetric riemannian two-manifold and ${\Bbb{G}}$ be the associated neutral K\"ahler metric 
on $TN$. Then the induced metric is degenerate on the graph $\xi\rightarrow (\xi,\eta=F(\xi,\bar{\xi}))$ with:
\[
F=\left[H(R)\pm i\left[B_2e^{-2u}+e^{-2u}\int\frac{(R\dot{H}-H)^2}{2R(1+R\dot{u})}e^{2u}dR\right]^{\scriptstyle{\frac{1}{2}}}\right]e^{i\theta}
\]
for any real differentiable function H.
\end{Prop}
\begin{pf}
The slopes $\lambda$ and $\sigma$ can be readily computed for this surface and it is then found that 
$\lambda^2=\sigma\bar{\sigma}$.
\end{pf}

\noindent{\bf Remark}: In presence of any Killing vector field on (N,g) we expect that 
a corresponding invariant area - stationary graph-type surface in TN exists and is given by a similar construction.

\section{The Space of Oriented Affine Lines in ${\Bbb{R}}^3$}

The four-manifold TS$^2$ can be identified with the space ${\Bbb{L}}$ of oriented affine lines in ${\Bbb{R}}^3$, and the
neutral K\"ahler metric ${\Bbb{G}}$, as constructed above with g equal to the round metric on S$^2$, is invariant under 
the action induced on ${\Bbb{L}}$ by the Euclidean isometry group acting on ${\Bbb{R}}^3$ \cite{gak4}. 

A surface $\Sigma$ in 
${\Bbb{L}}$ is a two-parameter family of oriented lines (or line congruence) in ${\Bbb{R}}^3$, which is the graph of a 
section of the bundle ${\Bbb{L}}\rightarrow$S$^2$ iff it can be parameterised by the direction $\xi$ of the oriented lines.
Moreover, a surface in ${\Bbb{L}}$ is lagrangian iff the associated line congruence is orthogonal
to a surface in ${\Bbb{R}}^3$. 

For the round metric $e^{2u}=4(1+\xi\bar{\xi})^{-2}$ and the above construction of area-stationary line congruences 
in ${\Bbb{L}}$ simplifies to:

\begin{Thm}
A rotationally symmetric surface in ${\Bbb{L}}$ which is given by the graph of a local section 
$\xi\rightarrow (\xi,\eta=F(\xi,\bar{\xi}))$ is area-stationary with respect to ${\Bbb{G}}$ iff
\[
F=\left[A_1R+B_1R^{-1}(1+R^2)^2\pm i\left[A_2R^2+B_2(1+R^2)^2-B_1^2R^{-2}(1+R^2)^4\right]^{\scriptstyle{\frac{1}{2}}}\right]e^{i\theta},
\]
for $A_1,A_2,B_1,B_2\in{\Bbb{R}}$, $A_2\neq0$, where $\xi=Re^{i\theta}$ and $g=e^{2u}d\xi d\bar{\xi}$ for $u=u(R)$.
\end{Thm}

To find closed area-stationary surfaces, we must have $B_1=0$, since otherwise the surface is asymptotic
to the fibre of the bundle ${\Bbb{L}}\rightarrow S^2$ at $R=0$. In addition, by a translation we can set $A_1=0$ and the
surface is determined by
\[
F=\pm i\left[B_2+C_2R^2+B_2R^4\right]^{\scriptstyle{\frac{1}{2}}}e^{i\theta},
\]
for $C_2\in{\Bbb{R}}$ such that $-2B_2\leq C_2$, and $B_2\geq0$.

This can be extended through $R=0$ and $R\rightarrow\infty$ and yields a two parameter family of area-stationary 
tori for $C_2\neq2B_2$. Under the projection map ${\Bbb{L}}\rightarrow$S$^2$ these tori double cover the sphere, 
except at the north and south pole, where the inverse image
of each of these points is a circle.

The induced metric is positive definite on the upper part and negative definite on the lower part of these 
tori, or vice versa, depending on the sign of $C_2-2B_2$. The inner and outer meridian circles (given by $R=1$) are 
totally null: the surface is both lagrangian and holomorphic at these points.

Finally for $C_2=2B_2$
\[
F=\pm i\sqrt{B_2}\left(1+R^2\right)e^{i\theta}
\]
is a torus on which the induced metric is degenerate.
\vspace{0.2in}

\noindent{\bf Acknowledgement}: The first author would like to thank Martin Stynes for drawing his attention to 
Lemma \ref{l:2ode}.

\end{document}